\documentclass[10pt]{amsproc}
\usepackage{amssymb,amsmath, amsthm, amsfonts}
\usepackage{mathrsfs}
\usepackage[all]{xypic}
\usepackage{graphicx}
\usepackage{latexsym}
\usepackage[mathscr]{euscript}
\usepackage{tikz,color} 
\usetikzlibrary{calc,matrix,arrows,chains,positioning,scopes,decorations.pathmorphing,backgrounds,fit}

\numberwithin{equation}{section}
\newtheorem{teo}{Theorem}[section]
\newtheorem{cor}{Corollary}[section]
\newtheorem{pro}{Proposition}[section]
\newtheorem{lem}{Lemma}[section]
\newtheorem{es}{Example}[section]
\newtheorem{defi}{Definition}[section]
\newtheorem{rem}{Remark}[section]
\newtheorem{notat}{Notation}[section]

\newcommand{\bdfn}{\begin{defi} \begin{rm}}
\newcommand{\edfn}{\end{rm} \end{defi}}
\newcommand{\bthm}{\begin{teo}}
\newcommand{\ethm}{\end{teo}}
\newcommand{\bprop}{\begin{pro}}
\newcommand{\eprop}{\end{pro}}
\newcommand{\bcor}{\begin{cor}}
\newcommand{\ecor}{\end{cor}}
\newcommand{\blem}{\begin{lem}}
\newcommand{\elem}{\end{lem}}
\newcommand{\bfact}{\begin{rem} \begin{rm}}
\newcommand{\efact}{\end{rm} \end{rem}}
\newcommand{\bex}{\begin{es} \begin{rm}}
\newcommand{\eex}{ \end{rm} \end{es}}
\newcommand{\bnot}{\begin{notat} \begin{rm}}
\newcommand{\enot}{\end{rm} \end{notat}}
\newcommand{\ten}{\otimes}
\newcommand{\quot}[2]{{\raisebox{.2em}{$#1$}\left/\raisebox{-.2em}{$#2$}\right.}}

\newcommand{\D}{\mathcal{D}}

\tikzset{node distance=3cm, auto}

\title{Notes on divisible MV-algebras}
\author{Serafina Lapenta\footnote{S. Lapenta: Department of Mathematics, University of Salerno, Via Giovanni Paolo II, 132 84084 Fisciano (SA), Italy. \textit{email:} slapenta@unisa.it},
Ioana Leu\c stean\footnote{I. Leu\c stean: Department of Computer Science, Faculty of Mathematics and Computer Science, University of Bucharest, Academiei nr.14, sector 1, C.P. 010014,  Bucharest, Romania, \textit{email:} ioana@fmi.unibuc.ro}
}

\begin{document}
\maketitle

\begin{abstract}
In these notes we study the class of divisible MV-algebras inside the algebraic hierarchy of MV-algebras with product. We connect  divisible MV-algebras with $\mathbb Q$-vector lattices, we present the divisible hull as a categorical adjunction and we prove a duality between finitely presented algebras and rational polyhedra. \\
\textit{Keywords:} DMV-algebras, MV-algebras, Rational \L ukasiewicz logic,  divisible hull, rational polyhedra.
\end{abstract}

\section*{Introduction}
MV-algebras were introduced by C.C. Chang \cite{Cha1} as algebraic counterpart of $[0,1]$-valued \L ukasiewicz logic. They are structures $(A, \oplus, ^*, 0)$ of type $(2,1,0)$ that satisfy appropriate axioms. 

Since their first appearance, MV-algebras have been deeply investigated and their theory has received a major boost after D. Mundici \cite{Mun1} proved a categorical equivalence with lattice-ordered Abelian groups with a strong order unit. 

In this work we will focus on the class of {\em divisible MV-algebras},  that correspond by the above mentioned categorical equivalence to divisible lattice-ordered groups with strong unit. 
Divisible MV-algebras are axiomatized and studied under the name of {\em DMV-algebras} by B. Gerla in \cite{DMV}.

Following similar results from the theory of lattice-ordered groups, any MV-algebra can be embedded in a divisble MV-algebra and any linearly ordered divisible MV-algebra is elementarily  equivalent to the standard MV-algebra on  $[0,1]$. These are  key steps for proving the completeness theorem for \L ukasiewicz logic \cite{Cha2,CDM}. 

 Different expansions of MV-algebras have been defined, following the natural hierarchy of lattice-ordered structures, from groups to algebras \cite{DiND,MonPMV,LeuRMV,LL1,LPhD}. In this paper we provide a slightly different, but equivalent, axiomatization  for DMV-algebras which allows us to include them in the general framework of MV-algebras with product. As consequence, in Section \ref{sec:ten}, we display an adjunction between semisimple MV-algebras and semisimple DMV-algebras via the  divisible hull functor, while in Section \ref{sec:logic} we connect DMV-algebras and $\mathbb Q$-vector lattices \cite{madden} with strong unit.
 
The propositional calculus that has DMV-algebras as models was defined in \cite{DMV}, under the name of {\em Rational \L ukasiewicz logic}. In Section \ref{sec:logic} we  provide an equivalent axiomatization  as a subsystem of the logic $\mathcal{R}$\L \ of Riesz MV-algebras, which are  MV-algebras closed for multiplication by scalars in $[0,1]$.

A recent line of research in the theory of MV-algebras is the investigation of geometrical dualities, following the steps of Baker and Beynon's duality theorem between finitely presented vector spaces and polyhedra \cite{Baker,Bey}. In \cite{LucaEnzo} one can find a duality theorem for finitely presented MV-algebras and rational polyhedra with $\mathbb Z$-maps, while in \cite{RMV-DiNola} finitely presented Riesz MV-algebras are proved to be dual to an appropriate category of polyhedra. In the final section of these notes we  complete the framework of geometrical dualities by proving that finitely presented DMV-algebras are dual to rational polyhedra with $\mathbb Q$-maps.

The results presented in this paper are expectable or considered "folklore" and some proofs are similar with already known ones.  However, our main goal is to present the divisible MV-algebras as MV-algebras endowed with a scalar multiplication with rational scalars from $[0,1]$ and to place them  inside the algebraic hierarchy of MV-algebras with product.  

\section{Preliminaries}\label{sec:prelim}
An MV-algebra is a structure $(A,\oplus,^*,0)$ which satisfies the following  properties for any $x,y \in A$:

(MV1) $(A,\oplus,0)$ is an Abelian monoid,

(MV2) $({x^*})^* = x$,

(MV3) $(x^* \oplus y)^* \oplus y = (y^* \oplus x)^* \oplus x,$

(MV4) $0^* \oplus x = 0^*$.\\
We refer to \cite{CDM} for all the unexplained notions related to MV-algebras.  In any MV-algebra $A$ we can define the following:
$1 = 0^*$,  $x \odot y = (x^* \oplus y^*)^*$,   
$x \vee y = (x \odot y^*) \oplus y$ and  $x \wedge y = (x \oplus y^*) \odot y$
for any $x,y \in A$. Hence $(A,\vee,\wedge, 0,1)$ is a bounded distributive lattice  such that $x\leq y$ if and only if $x\odot y^*=0$.\\

A Riesz MV-algebra \cite{LeuRMV} is a structure $(R, \oplus, ^*, 0, \{r\mid r \in [0,1]\})$ such that $(R, \oplus, ^*, 0)$ is an MV-algebra and $\{r\mid r \in [0,1]\}$ is a family of unary operators such that  the following  properties hold for any $x,y \in A$ and $r,q\in [0,1]$:

(RMV1) $r (x\odot y^{*})=(r  x)\odot(r y)^{*},$

(RMV2) $(r\odot q^{*})\cdot x=(r x)\odot(qx)^{*},$

(RMV3) $r (q  x)=(rq) x,$

(RMV4) $1 x=x.$\\

A PMV-algebra \cite{DiND,MonPMV} is a structure $(P, \oplus, \ ^*, \cdot , 0)$ such that $(P, \oplus, \ ^*, 0)$ is an MV-algebra, the operation  $\cdot:P\times P\to P $ is  associative and commutative,  and  the following identities hold for any $x,y,z \in P$:

(PMV1) $z\cdot (x\odot  y^*)=(z\cdot x)\odot (z\cdot y)^*$,

(PMV2) $x\cdot y=y\cdot x$,

(PMV3) $x\cdot 1=x.$\\
We remark that for this paper we are assuming any PMV-algebra to be commutative and unital, while the definition from \cite{DiND} is more general.\\

A DMV-algebra is an MV-algebra $A$ endowed with operators $\{\delta_n\}_{n\in \mathbb N}$ that satisfy the following, for any $x\in A$ and $n\in \mathbb{N}$:

(DMV1) $n\delta_nx=x$;

(DMV2) $\delta_n x \odot (n-1)\delta_n x =0$.\\
A \textit{divisible} MV-algebra is the MV-algebra reduct of a DMV-algebra.

\bfact \label{rem:div-DMV}
Note that is possible to provide an equivalent definition for divisible MV-algebras. In particular, an MV-algebra $A$ is divisible if for any $x\in A$ and for any $n\ge 0$ there exists a $y\in A$ such that $(n-1)y\odot y=0$ and $x=ny$. 
Any such algebra is indeed the reduct of a DMV-algebra, where $y=\delta_n(x)$.
\efact

Finally, Riesz MV-algebras are the particular case of the general construction of an MV-module over a PMV-algebra \cite{LeuMod}. In order to define this notion, we recall first the notion of \textit{partial sum} in an MV-algebra. We say that $x+y$ is defined if and only if $x\odot y=0$, and in this case $x+y=x\oplus y$. Then, if $P$ is a PMV-algebra and $A$ is a MV-algebra, $A$ is an MV-module over $P$ (or $P$-MV-module) if there is an external operation  $\varphi: P\times A\rightarrow A$, $\varphi(\alpha, x)=\alpha x $ such that, for any $x,y\in A$ and any $\alpha, \beta \in P$:

(MVM1) if $x+y$ is defined in $A$ then $\alpha x + \alpha y $ is also defined and $\alpha (x+y)= \alpha x+ \alpha y,$

(MVM2) if $\alpha + \beta $ is defined in $P$ then $\alpha x +\beta x $ is defined in \textit{A} and
$(\alpha + \beta )x=\alpha x + \beta x,$

(MVM3) $(\alpha \cdot \beta)x= \alpha (\beta x),$

(MVM4) $1_{P}x=x$ for any $x\in A$.\\
We remark that the definition of a MV-module has an equivalent equational form, hence the class of $P$-MV-modules forms a variety.\\

The standard model of an MV-algebra, PMV-algebra and Riesz MV-algebra is the unit interval $[0,1]$ endowed with the following operations: if $x, y,r\in [0,1]$ we define $x\oplus y=\min\{x+y,1\}$, $x^*=1-x$ and $x\cdot y=xy$ and $rx$ equal with the product of real numbers. Then $[0,1]$ naturally becomes an MV-algebra, a Riesz MV-algebra and a PMV-algebra, which generates the corresponding varieties in the first two cases. We remark that this is not longer true in the case of PMV-algebras. The standard model for a DMV-algebra is $[0,1]\cap\mathbb{Q}$, where $\delta_n(x)=\frac{x}{n}$ and it generates the variety of DMV-algebras.\\

A major boost in the theory of MV-algebras has been their categorical equivalence with Abelian  lattice-ordered groups with a strong unit. In the following we will provide some details on the aforementioned equivalence, and we urge the interested reader to consult \cite{CDM,DinLeu} for an in-depth treatment.

An {\em $\ell u$-group}  is a pair $(G,u)$, where $G$  is an Abelian lattice-ordered  group \cite{Birk} and $u$ is a strong unit. 
If $(G,u)$ is an $\ell u$-group, then $[0,u]_G=([0,u],\oplus, ^*, 0)$ is an MV-algebra, where  $[0,u]=\{ x\in G \mid 0\le x\le u\}$ and  $ x\oplus y=u\wedge (x+y)$, $x^*=u-x$ for any $x\in [0,u]$. 

If $\mathbf{MV}$ is the category of MV-algebras with MV-algebra homomorphisms and $\mathbf{uAG}$ is the category of $\ell u$-groups equipped with lattice-ordered group homomorphisms that preserve the strong unit, then one defines a functor $\Gamma: \mathbf{uAG}\rightarrow \mathbf{MV}$ by
$\Gamma (G,u)=[0,u]_G $ and $\Gamma(h)= h|_{[0,u_1]_{G_1}}$,
where  $(G,u)$ is an $\ell u$-group   and   $h: G_1\to G_2$ is a morphism in $\mathbf{uAG}$ between $(G_1,u_1)$ and $(G_2, u_2)$. In \cite{Mun1}, Mundici proved that  the functor $\Gamma$ establishes a categorical equivalence between  $\mathbf{uAG}$ and $\mathbf{MV}$.

All mentioned structures have an MV-algebra reduct. Hence, we can define forgetful functors from the categories $\mathbf{PMV}$ of PMV-algebras, $\mathbf{DMV}$ of DMV-algebras and $\mathbf{RMV}$ of Riesz MV-algebras to $\mathbf{MV}$. The categorical equivalence between MV-algebras and $\ell u$-groups can be generalized for each of this structures to an equivalence with an appropriate class of unital lattice-ordered structures having a lattice-ordered group reduct with a strong unit \cite{DiND,LeuRMV}.  In particular, in \cite{DMV}, DMV algebras are proved to be equivalent with divisible $\ell u$-groups.\\

In \cite{McN,LeuRMV,DMV} the free MV-algebra, Riesz MV-algebra and DMV-algebra are defined as algebras of $[0,1]$-valued functions. We will denote by $MV_X$, $RMV_X$ and $DMV_X$ the free algebra generated by $X$.  If $|X|=n$, the free algebras will be denoted by $MV_n$, $RMV_n$ and $DMV_n$ and they can be represented in terms of piecewise linear functions. In particular, $DMV_n$ is the algebra of functions in $n$-variables that are piecewise linear with rational coefficients.\\

An MV-algebra (DMV-algebra) is \textit{finitely generated} if it is generated by a finite set of elements, while it is \textit{finitely presented} if it is the quotient of a free finitely generated MV-algebra (DMV-algebra) by a finitely generated ideal. It is easy to see that for MV-algebras and DMV-algebras the notions of finitely generated ideal and principal ideal coincide. A \textit{projective} MV-algebra (DMV-algebra) is an algebra $A$ such that for any $B$ and $C$ MV-algebras (DMV-algebras) and any pair of homomorphisms $f:B \rightarrow C$, $g: A \rightarrow C$, with $f$ surjective, there exists $h: A\rightarrow B$ such that $g=f \circ h$.
It can be easily seen that free algebras are projective, see for example \cite[Example 17.2]{MunBook} for the case of MV-algebras.\\

Finally, an MV-algebra is semisimple if the intersection of its maximal ideals is $\{ 0 \}$. Equivalently, if it can be embedded in an algebra $C(X)$ of continuous functions $f:X\to [0,1]$, with $X$ compact Hausdorff space. Since a DMV-algebra has the same ideals of its MV-algebra reduct, a DMV-algebra is semisimple iff so is its MV-algebra reduct.

\section{DMV-algebras, $\mathbb{Q}$-vector lattices and their logic}\label{sec:logic}

We start this section giving an equivalent definition for a DMV-algebra. 
Upon setting $[0,1]\cap \mathbb{Q}=[0,1]_{\mathbb Q}$, we have the following preliminary result.

\bprop\label{pro:equivDefDMV}
Let $(D, \oplus, ^*, 0, \{r\mid r \in [0,1]_{\mathbb Q}\})$ be a structure such that $(D, \oplus, ^*, 0)$ is an MV-algebra and $\{r\mid r \in [0,1]_{\mathbb Q}\}$ is a family of unary operators. Then $D$ is an MV-module over $[0,1]_{\mathbb Q}$ if and only if the following properties hold for any $x,y \in A$ and $r,q\in [0,1]_{\mathbb Q}$:

\emph{(DMV1')} $r (x\odot y^{*})=(r  x)\odot(r y)^{*},$

\emph{(DMV2')} $(r\odot q^{*})\cdot x=(r x)\odot(qx)^{*},$

\emph{(DMV3')} $r (q  x)=(rq) x,$

\emph{(DMV4')} $1 x=x.$
\eprop
\begin{proof}
In \cite[Lemma 7 and Theorem 2]{LeuRMV} is proved that (RMV2) is equivalent with (MVM2) and (RMV1) is equivalent with (MVM1), when $r,q\in [0,1]$. In particular this holds for $r,q\in [0,1]_{\mathbb Q}$, that is, being (DMV1') and (DMV2') special cases of (RMV1) and (RMV2), the claim is settled. One can easily apply the techniques from \cite[Section 3]{LeuRMV} for a detailed proof.
\end{proof}

\bcor \label{cor:equivDefDMV}
Axioms (DMV1')-(DMV4') are equivalent with (DMV1)-(DMV2). That is, (DMV1')-(DMV4') is an equivalent definition for DMV-algebras.
\ecor
\begin{proof}
In \cite[Proposition 3.10]{LeuMod} is proved that $D$ is a divisible MV-algebra if and only if it is a MV-module over $[0,1]_{\mathbb Q}$. Hence, the conclusion follows from Remark \ref{rem:div-DMV} and Proposition \ref{pro:equivDefDMV}.
\end{proof}

\bfact
Following \cite[Proposition 2]{LeuRMV}, the previous characterization entails that any DMV-algebra must contain the interval $[0,1]_{\mathbb Q}$. We recall that an MV-algebra is simple iff it is a sub-algebra of $[0,1]$. Hence, any simple DMV-algebra is a dense subalgebra of $[0,1]$.
\efact

We now consider the following categories:

(i) $\mathbf{DMV}$, the category whose objects are DMV-algebras with homomorphism of DMV-algebras;

(ii) $\mathbf{uQVL}$, the category whose objects are $\mathbb{Q}$-vector lattice with a strong unit, i.e. lattice-ordered linear spaces over $\mathbb{Q}$, and whose maps are homomorphisms that preserve strong units, lattice order and scalar product.

\bthm 
The categorical equivalence defined by $\Gamma$ extends to $\mathbf{DMV}$and $\mathbf{uQVL}$. That is, the functor $\Gamma_\mathbb{Q}:\mathbf{DMV}\to \mathbf{uQVL}$ defined by $\Gamma_\mathbb{Q}(V,u)=[0,u]_V$, $\Gamma_\mathbb{Q}(f)=f\mid_{[0,u]_V}$, for any arrow $f: (V,u)\to (W,w)$ in $\mathbf{uQVL}$, establishes a categorical equivalence.
\ethm
\begin{proof}
It is straightforward from Corollary \ref{cor:equivDefDMV} and \cite[Theorem 4.6]{LeuMod},  where one can find the proof of a categorical equivalence between MV-modules over a PMV-algebra and lattice-ordered modules with a strong unit over a lattice-ordered ring with a strong unit. 
\end{proof}

Finally, we recall:
\bthm \cite[Theorem 4.5]{DMV} \label{teo:PWL}
The free DMV-algebra over $n$ generators, which is the algebra of term-functions from $[0,1]^n_Q\to [0,1]_{\mathbb Q}$, is isomorphic with the algebra of functions which are piecewise linear with rational coefficients.
\ethm

\bfact
The analogous of Theorem \ref{teo:PWL} for $\mathbb Q$-vector lattices is proved in \cite{madden}.
\efact

We end this section by describing the logic of DMV-algebras in terms of MV-modules over $[0,1]_{\mathbb Q}$. The \textit{Rational \L ukasiewicz logic} has been defined in \cite{DMV} as the logic that has DMV-algebras as models. We now provide a different logical system for DMV-algebras, denoted by $\mathcal{Q}$\L , and prove its equivalence with the one from \cite{DMV}.

We define the logic $\mathcal Q$\L \ as the one obtained by from \L ukasiewicz logic by adding a connective $\nabla_{r}$, for any $r\in [0,1]_{\mathbb Q}$. Axioms are the following.

(L1) $\varphi \rightarrow (\psi \rightarrow \varphi)$

(L2) $(\varphi \rightarrow \psi)\rightarrow((\psi \rightarrow \chi)\rightarrow(\varphi \rightarrow \chi))$

(L3) $ (\varphi \vee \psi)\rightarrow (\psi \vee \varphi)$

(L4) $(\neg \psi \rightarrow \neg \varphi)\rightarrow (\varphi \rightarrow \psi)$

(Q1) $\nabla_{r}(\varphi \rightarrow \psi)\leftrightarrow (\nabla_{r}\varphi \rightarrow \nabla_{r}\psi)$

(Q2) $ \nabla_{(r \odot q ^*)}\varphi \leftrightarrow (\nabla_{q}\varphi \rightarrow \nabla_{r}\varphi)$

(Q3) $\nabla_{r}(\nabla_{q}\varphi)\leftrightarrow \nabla_{r \cdot q} \varphi$

(Q4) $\nabla_1 \varphi \leftrightarrow \varphi$.\\
The only deduction rule is \textit{Modus Ponens}. Moreover, we define the connective $\Delta_{r}$ as $\neg\nabla_{r}\neg$, for any $r\in [0,1]$. We recall that Axioms (L1)-(L4) are the axioms of \L ukasiewicz logic.

We remark that for $r\in [0,1]$, the same axioms define the logical systems $\mathcal{R}$\L \ of Riesz MV-algebras \cite{LeuRMV}.

As usual, we can define the Lindenbaum-Tarski algebra of the logic. It is the quotient of the set of formulas $Form_{\mathcal{Q}}$ of the logic $\mathcal{Q}$\L \ by the equivalence relation defined as follows:
\begin{center}
$\varphi \equiv \psi$ if and only if $\vdash \varphi \rightarrow \psi$ and $\vdash \psi \rightarrow \varphi$
\end{center}
Let $Thm_{\mathcal{Q}}$ be the set of theorems of $\mathcal{Q}$\L . We define the following operations on $\quot{Form_{\mathcal{Q}}}{\equiv}$:
\begin{itemize}
\item[] $1=[Thm_{\mathcal{Q}}]$;
 $[\varphi ]^*=[\neg \varphi]$;
 $[\varphi]\oplus [\psi]=[\neg \varphi\rightarrow \psi]$;
\item[] $r [\varphi]=[\Delta_{r}\varphi]$ for $r\in {\mathbb Q} $.
\end{itemize}

\bthm
The structure $DMVL=(Form_{\mathcal{Q}}/\equiv, \oplus, *, {r|r\in [0,1]_{\mathbb Q}})$ is a DMV-algebra.
\ethm
\begin{proof}
It is a matter of computation to show that (Q1)-(Q4) are the logical expressions of the dual equations of (DMV1')-(DMV4'). A similar proof can be found in \cite[Proposition 5]{LeuRMV}.
\end{proof}

By the previous theorem, an appropriate semantics for the logic has DMV-algebras as models. If $A$ is a DMV-algebra, an evaluation is a function $e: Form_{\mathcal{Q}}\rightarrow A$ such that:

(e1) $e(\varphi \rightarrow \psi)= e(\varphi)^* \oplus e(\psi)$;

(e2) $e(\neg \varphi) =e(\varphi)^*$;

(e3) $e(\nabla_{r}\varphi)=(r e(\varphi)^*)^*$;\\
where $\oplus,\ *$ are the MV-algebra operations in $A$ and $r x$ the scalar product for any $r \in [0,1]_{\mathbb Q}$ and any $x\in A$.

Let $\varphi$ be a formula and $A$ a DMV-algebra. We say that $\varphi$ is a \textit{$A$-tautology} if $e(\varphi)=1$ for any evaluation $e:Form_{\mathcal{Q}}\rightarrow A$. We will denote a tautology by $\models_A \varphi$.

\bthm[Completeness]
Let $\varphi$ be a formula in $Form_{\mathcal{Q}}$, TFAE:

(1) $\vdash \varphi$;

(2) $\models _A \varphi$ for any DMV-algebra $A$;

(3) $\models _{[0,1]_{\mathbb Q}} \varphi$;

(4) $[\varphi]=1$ in $DMVL$.
\ethm
\begin{proof}
$(1)\Rightarrow (2)$ It is trivially seen that \textit{Modus Ponens} leads tautologies in tautologies. $(2)\Rightarrow (3)$ It is obvious. $(3)\Rightarrow (4)$ It follows from the fact that $[0,1]_{\mathbb Q}$ generates the variety of DMV-algebras \cite[Theorem 3.15]{DMV}. 
$(4)\Rightarrow (1)$ It follow directly by the definition of $1=[Thm_{\mathcal{Q}}]$.
\end{proof}

In \cite{DMV}, the \textit{Rational \L ukasiewicz logic} is defined as the logic that satisfies the natural translation of (DMV1)-(DMV2). That is, the Rational \L ukasiewicz logic is the expansion of \L ukasiewicz logic via a set of connectives $\{\delta_n \}_{n\in \mathbb{N}}$ satisfying the following axioms:

(L1),(L2), (L3), (L4)

(D1) $\stackrel{n}{\overbrace{\delta_n \varphi \oplus \dots \oplus \delta_n \varphi}} \rightarrow \varphi$

(D2) $\varphi \rightarrow \stackrel{n}{\overbrace{\delta_n \varphi \oplus \dots \oplus\delta_n \varphi}}$

(D3) $\neg \delta_n \varphi \oplus \neg (\stackrel{n-1}{\overbrace{\delta_n \varphi \oplus \dots \oplus\delta_n \varphi}})$.\\
We recall that the Rational \L ukasiewicz logic is complete w.r.t. $[0,1]_{\mathbb Q}$ \cite[Theorem 4.3]{DMV}.\\

We now define a correspondence between $\mathcal{Q}$\L \ and the Rational \L ukasiewicz logic. Let $ Form_{\mathcal{Q}}$ be the set of formulas in the former logic, and let $Form_{RatLuk}$ be the set of formulas in the latter.  To simplify our definition, we will use the connective $\Delta_r$, defined as $\neg \nabla_r \neg $, as primary.
\bdfn
(i) $\mathcal{I}_1: Form_{\mathcal{Q}}\to Form_{RatLuk}$, \\$\mathcal{I}_1(\delta_n \varphi)=\Delta_\frac{1}{n}\varphi$;

(ii) If $r=\frac{m}{n}$, $\mathcal{I}_2: Form_{RatLuk}\to Form_{\mathcal{Q}}$,\\
$\mathcal{I}_2(\Delta_r \varphi)=\stackrel{m}{\overbrace{\delta_n  \varphi \oplus \dots \oplus \delta_n \varphi} }$
\edfn 

\bthm
The correspondences $\mathcal{I}_1$ and $\mathcal{I}_2$ define two faithful translations of one logic into the other.
\ethm
\begin{proof}
We need to prove that both translations send theorems to theorems. Being both logics complete with respect to $[0,1]_{\mathbb Q}$, it is just a matter of computation to prove that (Q1)-(Q4) are send to tautologies in the Rational \L ukakasiewicz logic, and that (D1)-(D3) are send to tautologies in $\mathcal{Q}$\L .
\end{proof}

\section{DMV-algebras and the tensor product}\label{sec:ten}
In this section we display an adjunction between the categories of semisimple MV-algebras and semisimple DMV-algebras. In doing so, we provide a categorical framework for the construction of the divisible hull of a semisimple MV-algebra.

As in the previous section, we will denote by $[0,1]_{\mathbb Q}$ the interval $[0,1]\cap \mathbb{Q}$. 

We recall the definition of the semisimple MV-algebraic tensor product, as defined by D. Mundici \cite{Mun}. Let $A$, $B$ and $C$ be MV-algebras. A bimorphism is a map $\beta :A\times B \rightarrow C$ such that for $*\in \{ \vee , \wedge , + \}$,  $\beta (x, y_1 * y_2)=\beta (x,y_1) * \beta(x, y_2)$, $\beta (x_1 * x_2, y)=\beta (x_1,y)* \beta(x_2, y)$. Then, the tensor product $A\otimes_{MV} B$ is an MV-algebra and $\beta_{A,B}:A\times B\to A\otimes_{MV} B$ is its universal bimorphism such that for any MV-algebra $C$ and any bimorphism $\gamma: A\times B \rightarrow C$, there exists  an unique morphisms $\overline{\gamma}: A\ten_{MV}B \to C$ such that $\overline{\gamma}\circ\beta_{A,B}=\gamma$.

In \cite{Mun} the author proves that there exists a semisimple MV-algebra $A$ such that $A \ten_{MV}A$ is not semisimple. Therefore he defines the semisimple tensor product of $A$ and $B$, semisimple MV-algebras, by 
$$ A\ten B = \quot{A\ten_{MV}B }{ Rad(A\ten_{MV}B)}.$$ 
\noindent Consequently, $a\ten b$ is the class of $a\ten_{MV} b$ for any $a\in A$ and $b\in B$. For the semisimple tensor product a universal property holds  with respect to semisimple  MV-algebras, namely the following:

{\em for any semisimple MV-algebra $C$ and for any bimorphism $\beta: A\times B\rightarrow C$, there is a unique homomorphism of MV-algebras $\omega :A\otimes B\rightarrow C$ such that $\omega \circ \beta_{A,B}=\beta$.}\\

The following lemma ensures that a DMV-algebra has the same homomorphisms of its MV-algebra reduct.
\blem\label{lem:homom}
If $A$, $B$ are DMV-algebras and $f:A \to B$ is an MV-algebra homomorphism between the underlying MV-algebra reducts, then $f$ is a homomorphism of DMV-algebras.
\elem
\begin{proof}
We first remark that condition (ii) in the definition of a DMV-algebra implies that, for any $n$, we get  $\stackrel{n}{\overbrace{\delta_n x \oplus \dots \oplus \delta_n x}}\ =\ \stackrel{n}{\overbrace{\delta_n x + \dots + \delta_n x}}$, i.e. the MV-sum $\oplus$ coincides with the partial sum $+$. Then, we have the following:
\vspace*{-0.2cm}
\begin{center}
$n\delta_n(f(x))=f(x)=f(n\delta_n(x))=nf(\delta_nx)$.
\end{center}
By the above remark, the claim follows from cancellation of $+$, see \cite[Lemma 1.1(i)]{LeuMod}.
\vspace*{-0.2cm}
\end{proof}

It is well known that any MV-algebra can be embedded in a divisible one (see, for example, \cite[Corollary 5.3.1]{DinLeu}). The following proposition, together with \cite[Proposition 2.1]{LLTP1}, gives a different way to construct the embedding for semisimple MV-algebras.
\bprop\label{pro:tenDMV}
If $A$ is a semisimple MV-algebra,  then $[0,1]_{\mathbb Q} \ten A$ is a semisimple DMV-algebra.
\eprop
\begin{proof}
\cite[Theorem 3.1]{LLTP1} and \cite[Proposition 3.10]{LeuMod} entail that $[0,1]_{\mathbb Q} \ten A$ is a semisimple divisible MV-algebra. By Remark \ref{rem:div-DMV} and the fact that ideals in a DMV-algebras are MV-ideals, $[0,1]_{\mathbb Q} \ten A$ is a semisimple DMV-algebra.
\end{proof}

Let $\mathcal{U}_{\delta}$ be the forgetful functor from the full subcategory $\mathbf{DMV_{ss}}$ of semisimple DMV-algebras to the full subcategory $\mathbf{MV_{ss}}$ of semisimple MV-algebras.

\bprop \label{pro:maps}
Let $B$ be a  semisimple MV-algebra. For any  semisimple DMV-algebra $V$ and for any homomorphism of MV-algebras $f:B\rightarrow \mathcal{U}_{\delta}(V)$ there is a unique homomorphism of DMV-algebras
$\widetilde{f}:[0,1]_{\mathbb Q}\ten B\rightarrow V$ such that $\widetilde{f}\circ \iota_{B}=f$.
\eprop
\begin{proof}
Define $\beta_f:[0,1]_{\mathbb Q}\times B\to V$ by $\beta_f(\alpha,x)=\alpha f(x)$ for any $\alpha\in [0,1]_{\mathbb Q}$ and $x\in B$. Since $[0,1]_{\mathbb Q}$ is totally ordered and $V$ is a DMV-algebra (hence, an MV-module over $[0,1]_{\mathbb Q}$, by Corollary \ref{cor:equivDefDMV}), it is easily seen that $\beta_f$ is a bimorphism. Indeed we have

i) $\beta_f(\alpha_1 + \alpha_2,x)=(\alpha_1 + \alpha_2) f(x)=\alpha_1 f(x)+\alpha_2 f(x)=\beta_f(\alpha_1,x)+\beta_f(\alpha_2,x)$;

ii) $\beta_f(\alpha_1 \wedge \alpha_2,x)=(\alpha_1 \wedge \alpha_2) f(x)$. Since $[0,1]_{\mathbb Q}$ is totally ordered, without loss of generality, let $\alpha_1 \wedge \alpha_2=\alpha_1$. Then $(\alpha_1 \wedge \alpha_2) f(x)=\alpha_1 f(x)= \alpha_1 f(x) \wedge \alpha_2 f(x)$ by monotonicity of scalar multiplication. Being the latter equal to $\beta_f(\alpha_1, x)\wedge \beta_f(\alpha_2, x)$ we get the desired conclusion, which follows analogously for $\vee$.

iii) The proof that $\beta_f(\alpha, x_1+x_2)=\beta_f(\alpha,x_1)+\beta_f(\alpha, x_2)$, $\beta_f(\alpha, x_1\wedge x_2)=\beta_f(\alpha,x_1)\wedge \beta_f(\alpha, x_2)$ and $\beta_f(\alpha, x_1\vee x_2)=\beta_f(\alpha,x_1)\vee \beta_f(\alpha, x_2)$ is similar to (i), via \cite[Definition 3.1 and Lemma 3.19]{LeuMod}.\\
We now use the universal property of the semisimple tensor product.  We note that, by Lemma  \ref{lem:homom}, any MV-homomorphism between DMV-algebras is an homomorphism of DMV-algebras.
\end{proof}

We are now ready to define a functor ${\mathcal D}_\ten:\mathbf{MV_{ss}}\to \mathbf{DMV_{ss}}$ by 

${\mathcal D}_\ten(B)=[0,1]_{\mathbb Q}\ten B$ for any semsimple MV-algebra $B$ and 

${\mathcal D}_\ten(f) =\widetilde{f}$ for any homomorphism of MV-algebras $f:A\to B$, where $\widetilde{f}:[0,1]_{\mathbb Q}\ten A\to [0,1]_{\mathbb Q}\ten B$ is the unique homomorphism of DMV-algebras such that $\widetilde{f}\circ \iota_A=\iota_B\circ f$, which exists by Proposition \ref{pro:maps}.

\bthm\label{teo:adjMV-DMV}
Under the above hypothesis,  ($\mathcal{D}_\ten, \mathcal{U}_{\delta})$  is an adjoint pair. 
\ethm
\begin{proof}
The proof follows \textit{mutati mutandis} from \cite[Theorem 4.2]{LLTP1}.
\end{proof}

\bfact
In a similar fashion, we can display an adjunction between $\mathbf{DMV_{ss}}$ and $\mathbf{RMV_{ss}}$, the full subcategory of semisimple Riesz MV-algebras. We consider the following functors:

(i) $\mathcal{T}_{\ten}^{\mathbb{Q}} : \mathbf{DMV_{ss}} \to \mathbf{RMV_{ss}}$ defined as $\mathcal{T}_{\ten}^{\mathbb{Q}}(A)=A\ten[0,1]$ and on maps by the analogous of Proposition \ref{pro:maps};

(ii) $\mathcal{U}_\mathbb{Q}:\mathbf{DMV_{ss}}\to \mathbf{DMV_{ss}} $ is the forgetful functors that consider only multiplication with rational numbers.\\
One can easily prove that this is again a pair of adjoint functors, via the results proved in \cite{LLTP1} for Riesz MV-algebras.
\efact

In \cite[Remark 2.4 and Lemma 2.2]{LeuDia} the \textit{divisible hull} $A^d$ of a semisimple MV-algebra $A$ is defined. If $A\subseteq C(X)$, for a suitable compact Hausdorff space, then $A^d=\{a\in C(X)\mid a= \frac{a_1}{n}+ \dots + \frac{a_n}{n} \text{ for some }n\in \mathbb{N} \text{ and }a_1, \dots , a_n \in A  \}$. Let us denote by $\iota_d$ the embedding of $A$ in its divisible hull $A^d$.

\blem \label{univ-prop-divhull}
For any MV-algebra $A$, the following hold.

(i) $A^d=<\iota_d(A)>_{DMV}$,

(ii) The embedding $\iota_d: A \to A^d$ is essential,

(iii) For any DMV-algebra $D$ and any $f:A \to \mathcal{U}_{\delta}(D)$ there exists a unique $f_d:A^d\to D$ such that $f_d\circ \iota_d=f$. If $f$ is an embedding, $f_d$ is an embedding.
\elem
\begin{proof}
(i) Upon identifying $\iota_d(A)$ with $A$,  any $a\in A^d$ can be written in $<A>_{DMV}$ as $\delta_n(a_1) \oplus \dots \oplus \delta_n(a_n)$ for an appropriate $n\ge 1$. On the other side, an element in $<A>_{DMV}$ is inductively construct from elements of $A$ via $\delta_n$. If $a=\delta_n(b)$, then $a=\frac{b}{n}+ \frac{0}{n}\dots + \frac{0}{n}$ and it belongs to $A^d$ by definition.

(ii) and (iii) are proved in \cite[Lemma 2.2]{LeuDia}.
\end{proof}

\bprop \label{pro:divhull}
If $A$ is a semisimple MV-algebra and $D$ is a semisimple DMV-algebra such that 
$A\subseteq D$ and $D=\langle A\rangle_{DMV}$, then $D\simeq A^d$. 
\eprop
\begin{proof}
By Lemma \ref{univ-prop-divhull} there exists $e:A^d\to D$ such that $e(\iota_d(a))=a$ for any $a\in A$. Then we have
\begin{center}
$e(A^d)=e(<\iota_d(A)>_{DMV})=<e(\iota_d(A))>_{DMV}=$\\$=<A>_{DMV}=D$,
\end{center}
and $e$ is an isomorphism.
\end{proof}
\bcor
If $A$ is a DMV-algebra, then $A\simeq A^d$.
\ecor

We now recall the following theorem.
\bthm \cite[Theorem 4.2.4]{LPhD}\label{teo:univTen}
Let $X$ be a set, $<X>_{MV}$ the MV-algebra generated by $X$ and $P$ a totally ordered PMV-algebra. The following properties hold:\\
 (a)  $P \ten <X>_{MV}$ is generated as a $P$-MV-module by $\{1\ten x:x\in X\}$,
 
 (b)  if $M$ is a $P$-MV-module then
 for any function $f:X\rightarrow M$ there is a unique homomorphism of $P$-MV-modules
 $\widetilde{f}:P \ten MV_X \rightarrow M$ such that
 \begin{center}
  $\widetilde{f}(1\ten x)=f(x)$ for any $x\in X$.
\end{center}
  \ethm

Being $[0,1]_{\mathbb Q}$ a totally ordered PMV-algebra, a direct consequence of Theorem \ref{teo:univTen} is the following corollary.
\bcor
For any semisimple MV-algebra $A$,\\ $[0,1]_{\mathbb Q} \ten A$ is the divisible hull of $A$.
\ecor
\begin{proof}
It follows from Proposition \ref{pro:divhull}, Theorem \ref{teo:univTen} and the remark that $\{1\ten x:x\in X\}=\iota_A(A)$, where $\iota_A:A\to [0,1]_{\mathbb Q}\ten A$, defined by $a\mapsto 1\ten a$, is the embedding in the MV-algebraic tensor product.
\end{proof}

Finally, we have the following.
\bthm \label{teo:freeTensDMV}
For any nonempty set $X$, $[0,1]_{\mathbb Q}\ten MV_X \simeq (MV_X)^d\simeq DMV_X$
\ethm
\begin{proof}
Assume $D$ is a semisimple DMV-algebra and  $f:X\to D$ is a function. Hence there is a unique homomorphism of MV-algebras 
$\overline{f}: MV_X\to {\mathcal{ U}}_{\delta}(D)$  which extends $f$.  By Proposition \ref{pro:maps}, there exists a homomorphism of DMV-algebras $\widetilde{f}:[0,1]_{\mathbb Q}\ten MV_X \to D$ such that $\widetilde{f}\circ\iota_{MV_X}=\overline{f}$, so $\widetilde{f}(1\ten x)=f(x)$ for any $x\in X$. The uniqueness of $\widetilde{f}$ is a consequence of the uniqueness of $\overline{f}$. Since $\iota_{MV_X}$ is an embedding we have $X\simeq \{1\ten x\mid x\in X\}$ and $[0,1]_{\mathbb Q}\ten MV_X$ is the free object in $\mathbf{DMV_{ss}}$. Being the free DMV-algebra $DMV_X$ an object in $\mathbf{DMV_{ss}}$, $[0,1]_{\mathbb Q}\ten MV_X \simeq DMV_X$.
\end{proof}

By categorical equivalence, the adjunctions 
$(\mathcal{T}_{\ten}^{\mathbb{Q}} , \mathcal{U}_{\mathbb{\mathbb{Q}}})$ and $(\mathcal{D}_\ten , \mathcal{U}_{\delta})$ can be naturally transferred to lattice-ordered structures. We denote by $\mathbf{auG_a}$ the category of archimedean $\ell u$-groups; $\mathbf{uQVL_a}$ the category of archimedean $\mathbb{Q}$-vector lattices with a strong unit; $\mathbf{uRS_a}$ the category of archimedean Riesz Spaces with strong unit.

In \cite{LLTP1}, the adjunction $(\mathcal{T}_\ten , \mathcal{U}_{\mathbb{R}})$ between semisimple MV-algebras and semisimple Riesz MV-algebras with a strong unit is proved, and it is extended to $(\mathcal{T}_{\ten a} , \mathcal{U}_{\ell \mathbb{R}})$, the latter being an adjunction between $\mathbf{auG_a}$ and $\mathbf{uRS_a}$. 

With the same ideas, applying the inverse of  $\Gamma$ and $\Gamma_{\mathbb Q}$ (as defined in \cite{DMV} for MV-algebras), $(\mathcal{D}_\ten , \mathcal{U}_{\delta})$ extends to $(\mathcal{D}_{\ten a} , \mathcal{U}_{\ell \delta})$. This is an adjunction between $\mathbf{auG_a}$ and $\mathbf{uQVL_a}$.

Applying the converses of the functors $\Gamma_{\mathbb{Q}}$ and $\Gamma_{\mathbb{R}}$ , $(\mathcal{T}_{\ten}^{\mathbb{Q}} , \mathcal{U}_{\mathbb{Q}})$ extends to $(\mathcal{T}_{\ten a}^\mathbb{Q} , \mathcal{U}_{\mathbb{Q} a})$. This is an adjunction between $\mathbf{uQVL_a}$ and $\mathbf{uRS_a}$.

The following diagram summarizes our results.
\begin{center}
\begin{tikzpicture}
  \node (A) {$\mathbf{auG_a}$};
  \node (B) [below of=A] {$\mathbf{MV_{ss}}$};
  \node (C) [right of=A] {$\mathbf{uQVL_{a}}$};
  \node (D) [below of=C] {$\mathbf{DMV_{ss}}$};
  \node (E) [right of=C] {$\mathbf{uRS_a}$};
  \node (F) [below of=E] {$\mathbf{RMV_{ss}}$};
  \draw[->] (A) to node [swap] {$\Gamma$} (B);
  \draw[->] (C) to node [swap] {$\Gamma_\mathbb{Q}$} (D);
  \draw[->] (B) to node  {$\mathcal{D}_\ten$} (D);
  \draw[->] (E) to node {$\Gamma_{\mathbb{R}}$} (F);
  \draw[->] (D) to node {$\mathcal{T}_{\ten}^{\mathbb{Q}}$} (F);
  \draw[->, bend right] (B) to node  {$\mathcal{T}_{\ten}$} (F);
  \draw[->, bend left] (A) to node [swap] {$\mathcal{T}_{\ten a}$} (E);
  \draw[->] (C) to node [swap]{$\mathcal{T}_{\ten  a}^{\mathbb{Q}}$} (E);
  \draw[->] (A) to node [swap] {$\mathcal{D}_{\ten a}$} (C);

\end{tikzpicture}
\end{center}

\section{DMV-algebras and rational polyhedra}\label{sec:dualDMV}
We start this section recalling all needed notions of polyhedral geometry.

A $m$-simplex in $[0,1]^n$ is the convex hull $C$ of $m+1$ affinely independent points $\{ v_o, \ldots , v_m\}$ in the euclidean space $[0,1]^n$; the points $v_i$ are called vertexes and $C$ is rational if any vertex has rational coordinates.  A (\textit{rational}) \textit{polyhedron} is the union of finitely many (rational) simplexes; any simplex is obtainable as a finite union of a finite intersection of closed half-spaces.\\
If $P \subseteq [0,1]^n$ is a polyhedron, a $\mathbb{Z}$-map $z: P\rightarrow [0,1]^m $ is a continuous map $z=(z_1, \ldots z_m)$ where any $z_i:[0,1]^n\rightarrow [0,1]$ is a piecewise linear function with integer coefficient.

Rational polyhedra with $\mathbb{Z}$-maps form a category which is dually equivalent (see \cite{LucaEnzo}) to the full subcategory of finitely presented MV-algebras. Moreover, finitely presented projective MV-algebras are in duality with retracts \--- by $\mathbb{Z}$-maps \--- of finite-dimensional unit cubes. This duality follows from the ideas of Baker and Beynon \cite{Baker,Bey}. We now prove the analogous duality theorem. Most of the results of this section can be proved exacly as in \cite{CDM,RMV-DiNola,LucaEnzo}, relying on subsequently remark.

\bfact \label{rem:idRMV}
Let $I$ be a principal ideal in $DMV_n$, generated by $a\in DMV_n\setminus MV_n$. By Theorem \ref{teo:freeTensDMV} there exists a natural number $k$ and $a_1,\ldots, a_k\in MV_n$ such that $a=\frac{1}{k}a_1\oplus\cdots \oplus\frac{1}{k}a_k$. Let $b=a_1\oplus\cdots \oplus a_k$. Since $\frac{1}{k}a_i \le a_i$, we have $a\le b$ therefore the ideal generated by $a$ is included in the ideal generated by $b$, in symbols $(a] \subseteq (b]$. For the other inclusion, $\delta_k(a_i)=\frac{1}{k}a_i\le a$ for any $i$, therefore $a_i=k\delta_k(a_i)\in (a]$. Hence $b\in (a]$,  and $(a] = (b]$. 

This means that we can replace the generator of a principal ideal in $DMV_n$ by an element of $MV_n$.
\efact

In analogy with the case of MV-algebras and Riesz MV-algebras, we have the following definitions.
\begin{itemize}
\item[(i)] Given subset $S\subseteq DMV_n$,  $V(S)=\{ \mathbf{x}\in [0,1]^n\mid f(\mathbf{x})=0 \mbox{ for any }f\in S\}$.
\item[(ii)] Given a subset $X\subseteq [0,1]^n$,  $I(X)=\{ f\in DMV_n\mid f(\mathbf{x})=0 \mbox{ for any } \mathbf{x}\in X\}$.
\end{itemize}
Since functions in $DMV_n$ are continuous, $f^{-1}(0)$ is a closed set and being $V(S)=\bigcap _{f\in S}f^{-1}(0)$, any $V(S)$ is closed in $[0,1]^n$.

We recall that a subalgebra $A$ of $C([0,1]^n)$ is called \textit{separating} if for any two point $\mathbf{x}, \mathbf{y}\in [0,1]^n$ there exists $f\in A$ such that $f(\mathbf{x})=0$ and $f(\mathbf{y})>0$.
\blem\label{lem:separating}
$DMV_n$ is a separating algebra of $C([0,1]^n)$.
\elem
\begin{proof}
In \cite[Lemma 3.4.6]{CDM} is proved that $MV_n$ is separating, for any $n\in \mathbb{N}$. The claim follows from the remark that $MV_n\subseteq DMV_n$.
\end{proof}

The following proposition is a direct consequence of the analogous results in \cite[Chapter 3]{CDM}, Lemma \ref{lem:separating} and the fact that homomorphisms and ideals of a DMV-algebra are the same of its MV-algebra reduct.
\bprop \label{pro:CDM}
The following properties hold:

(i)For each closed subset $C$ of $[0,1]^n$, $C=V(I(C))$.

(ii) Let $C$ be a closed subset of $[0,1]^n$. For any $\mathbf{x}\in C$ the map $\mathbf{x}\to I(\mathbf{x})$ is an isomorphism between $C$ and $Max(DMV_n\mid_C)$

(iii) Let $J$ be a proper ideal. Then $I(V(J))$ is the intersection of the maximal ideals containing $J$.

(iv) An ideal $J$ in $DMV_n$ is an intersection of maximal ideals if and only if $J=I(V(J))$.
\eprop

\blem \label{lem:inclID}
Let $f,g \in DMV_n$. Then, $g\in (f]$ if and only if $V(f)\subseteq V(g)$
\elem
\begin{proof}
One direction is trivial. For the other direction, let $V(f)\subseteq V(g)$. By Remark \ref{rem:idRMV} there exists $g^*$ and $f^*$ in $MV_n$ such that $(f]=(f^*]$ and $(g]=(g^*]$. Moreover, by construction of $f^*$, $f(x)=0$ iff $f^*(x)=0$, that is $V(f)=V(f^*)$ and $V(g)=V(g^*)$, hence $V(f^*)\subseteq V(g^*)$. Therefore, by \cite[Lemma 3.4.8]{CDM}, $g^*$ belongs to $(f^*]_{MV}$, the ideal generated by $f^*$ in $MV_n$ and trivially, it belongs to the ideal it generates in $DMV_n$. Hence, $g\in (g^*] \subseteq  (f^*]=(f]$.
\end{proof}

\bprop \label{pro:maxInt}
Any principal ideal in $DMV_n$ is intersection of maximal ideals.
\eprop
\begin{proof}
Let $J=(f]_{DMV}$ be a principal ideal in $DMV_n$. We first remark that $V(J)=V(f)$. Since the inclusion $J\subseteq I(V(J))$ always holds, let us prove the converse inclusion.\\
Let $g \in I(V(J))$. If $x\in V(J)$, we have $g(x)=0$. Thus $x\in V(g)$ and $V(J)=V(f)\subseteq V(g)$. By Lemma \ref{lem:inclID}, $g\in J$ and the conclusion follows from Proposition \ref{pro:CDM}(iv).
\end{proof}

\blem \label{semisimple}
A finitely generated DMV-algebra $D=\quot{DMV_n}{J}$ is semisimple if and only if $J$ is intersection of maximal ideals of $DMV_n$.
\elem
\begin{proof}
We consider the projection 
\vspace*{-0.2cm}
\begin{center}
$\pi_D: DMV_n \to \quot{DMV_n}{J}.$
\end{center}
\vspace*{-0.2cm}
Let us assume that $D$ is semisimple. An ideal is maximal in $D$ if it is of the form $\quot{M}{J}$, with $M$ maximal ideal of $DMV_n$ that contains $J$ \-- the proof can be easily deduced from \cite[Proposition 1.2.10]{CDM} \-- since $J$ is semisimple, $ \bigcap_{\quot{M}{J}\in Max(D)} \quot{M}{J}=\{0 \}$. Then, 
\begin{eqnarray}
J&=&\pi_D^{-1}(0)=\pi_D^{-1}\left( \bigcap_{\quot{M}{J}\in Max(D)} \quot{M}{J}\right)=\nonumber \\
&=&\bigcap_{M\in Max(DMV_n), J\subseteq M} \pi_D^{-1}(\quot{M}{J})=\nonumber\\
&=&\bigcap_{M\in Max(DMV_n), J\subseteq M} M.\nonumber
\end{eqnarray}
On the other direction, we can assume that $J$ is the intersection of all maximal ideals that contain it. Denoted by $\{ M_i \}_{i\in I}$ the collection  of such ideals, we get $Max\left(\quot{DMV_n}{J}\right)=\{\pi_J(M_i)\}_{i\in I}$. Hence $\bigcap_i \pi_J(M_i)=\pi_J(J)=0$ and $D$ is semisimple.
\end{proof}

\bprop
Any finitely presented DMV-algebra is semisimple.
\eprop
\begin{proof}
It follows from Proposition \ref{pro:maxInt} and Lemma \ref{semisimple}.
\end{proof}

\bprop \label{pro:CDM2}
For any closed subset $C$ of $[0,1]^n$, $DMV_n\mid_C$ is isomorphic with $\quot{DMV_n}{I(C)}$.
\eprop
\begin{proof}
It follows from Lemma \ref{lem:separating}, Proposition \ref{pro:CDM} and \cite[Proposition 3.4.5]{CDM}.

\end{proof}

The following theorem is crucial in the proof of duality.
\bthm \label{thm:zerosetsDMV}
The zerosets of elements of $DMV_n$ are exactly rational polyhedra.
\ethm
\begin{proof}
Let $P$ be a rational polyhedron. By \cite[Lemma 3.2]{LucaEnzo} there exists $f\in MV_n\subseteq DMV_n$ such that $P=f^{-1}(0)$.\\
On the other direction, let $f$ be an element of $DMV_n$. By Remark \ref{rem:idRMV} there exists $g\in MV_n$ such that $(f]_{DMV}=(g]_{DMV}$. By Lemma \ref{lem:inclID} it follows that $f^{-1}(0)=V(f)=V(g)=g^{-1}(0)$. This is a rational polyhedron by \cite[Lemma 3.2]{LucaEnzo}.
\end{proof}

\bprop \label{prop:reprDMV}
Let $J$ be a principal ideal in $DMV_n$. Then there exists a rational polyhedron $P\subseteq [0,1]^n$ such that \  $\quot{DMV_n}{J}\simeq DMV_n\mid_P$, the set of piecewise linear functions with rational coefficients restricted to $P$.
\eprop
\begin{proof}
Let $f$ be the generator of $J$. By Theorem \ref{thm:zerosetsDMV}, we set $P=V(f)=V(J)$. By Proposition \ref{pro:maxInt} and Proposition \ref{pro:CDM}, $I(P)=I(V(J))=J$ and the claim is settled by Proposition \ref{pro:CDM2}.
\end{proof}

In the sequel, let  $\mathbf{DMV_{fp}}$ be the full subcategory of finitely presented DMV-algebras and homomorphisms of DMV-algebras, and let $\mathbf{RatPol_{[0,1]}^{\mathbb Q}}$ be the category of rational polyhedra laying in some appropriate unit cube with unital $\mathbb{Q}$-maps, that is maps $\lambda=(\lambda_1, \dots , \lambda_m):P \subseteq [0,1]^n\to Q\subseteq [0,1]^m$ in which the components $\lambda_i$ are elements of the free DMV-algebra $DMV_n$, i.e. continuous piecewise linear functions from $[0,1]^n$ to $[0,1]$ with rational coefficients. More precisely, any $\lambda$ is the restriction to $P$ and the co-restriction to $Q$ of the map $\lambda:[0,1]^n\to [0,1]^m$.

It is straightforward that $\mathbf{RatPol_{[0,1]}^{\mathbb Q}}$ is a category, since the image of a rational polyhedron via a piecewise linear function with rational coefficient is again a rational polyhedron. Moreover, we have the following.

\blem
Let $f:[0,1]^n \to [0,1]$ be a piecewise linear function and let $\{ f_1, \dots f_k \}$ its components, with $f_i(\mathbf{x})=a_{i_1}x_1 + \dots a_{i_n}x_n$ for $i=1, \dots , k$ and $\mathbf{x} \in [0,1]^n$. 

Let $\lambda=(\lambda_1, \dots , \lambda_n):[0,1]^m \to [0,1]^n$ be a $\mathbb{Q}$-map and for each $\lambda_l$, let $\{ g_{l_1}, \dots g_{l_{t_l}} \}$ the set of its components, where $g_{l_j}(\mathbf{x})=b_{j_1}x_1+ \dots + b_{j_m}x_m$, $\mathbf{x}\in [0,1]^m$.

Then $f\circ \lambda$ is a piecewise linear function. Moreover, the composition of unital $\mathbb{Q}$-maps is a unital $\mathbb{Q}$-map.
\elem
\begin{proof}
For any $\mathbf{x}\in [0,1]^m$, there exist $f_i, g_{p_{s_p}}$ for $p=1, \dots, n$ such that $(f\circ \lambda )(x)=f_i(g_{1_{s_1}}(\mathbf{x}), \dots, g_{n_{s_n}}(\mathbf{x})) $  for an appropriate choice of the components. Hence, 
\begin{eqnarray}
&&f_i(g_{1_{s_1}}(\mathbf{x}), \dots, g_{n_{s_n}}(\mathbf{x}))=\nonumber \\
&=& a_{i_1}g_{1_{s_1}}(\mathbf{x})+ \dots + a_{i_n}g_{n_{s_n}}(\mathbf{x})=\nonumber\\
&=&a_{i_1}(\sum_{h=1}^m b_{s_{1_h}}x_h)+ \dots + a_{i_n}((\sum_{h=1}^m b_{s_{n_h}}x_h)).	\nonumber
\end{eqnarray}
Re-arranging each term, we obtain an affine function with rational coefficients.

Finally, if $\lambda: P\subseteq [0,1]^n\to Q\subseteq [0,1]^m$ and $\sigma: Q\subseteq [0,1]^m\to R\subseteq [0,1]^s$, we have that 
\begin{eqnarray}
&&(\sigma \circ \lambda)(\mathbf{x})=(\sigma (\lambda_1(\mathbf{x}), \dots , \lambda_m(\mathbf{x})))= \nonumber\\
&=&(\sigma_1 (\lambda_1(\mathbf{x}), \dots , \lambda_m(\mathbf{x})), \dots , \sigma_s(\lambda_1(\mathbf{x}), \dots , \lambda_m(\mathbf{x})))= \nonumber \\
&=& ((\sigma_1\circ \lambda)(\mathbf{x}), \dots (\sigma_s\circ \lambda)(\mathbf{x})). \nonumber
\end{eqnarray}
By the previous step, $(\sigma_i\circ \lambda)$ is a piecewise linear function with rational coefficients and $(\sigma \circ \lambda)$ is a $\mathbb{Q}$-map.
\end{proof}

\noindent We define a functor $\D :\mathbf{RatPol_{[0,1]}^{\mathbb Q}}\rightarrow \mathbf{DMV_{fp}}$ by
\begin{itemize}
\item[] $\D (P)=DMV_n\mid_{P}$, for any polyhedron $P\subseteq [0,1]^n$;
\item[] $\D (\lambda): \D (Q) \rightarrow \D (P)$ defined by $\D (\lambda)(f)=f \circ \lambda$, for any $\lambda: P \rightarrow Q$, with $P\subseteq [0,1]^n$ and $Q\subseteq [0,1]^m$.
\end{itemize}

\noindent It is easily seen (cfr \cite[Lemma 3.3]{LucaEnzo}) that $\D (\lambda)$ is an homomorphism of MV-algebras and by Lemma \ref{lem:homom}, it is an homomorphism of DMV-algebras. Trivially, $\D$ is a functor.

\bthm \label{teo:BBfpDMV}
The functor $\D$ establishes a categorical equivalence between the category  $\mathbf{RatPol_{[0,1]}^{\mathbb Q}}$ and the opposite of the category $\mathbf{DMV_{fp}}$.
\ethm
\begin{proof}
The proof follows \textit{mutati mutandis} from \cite[Theorem 3.4]{LucaEnzo}. We include all details for sake of completeness.

We need to prove that $\D$ is full, faithful and essentially surjective. Proposition \ref{prop:reprDMV} entails that $\D$ is essentially surjective.

To show that $\D$ is faithful, let $\lambda, \mu : P\subseteq [0,1]^n\rightarrow Q\subseteq [0,1]^m$ such that $\lambda \neq \mu$. Therefore there exists $p\in P$ such that $(\lambda_1(p), \dots , \lambda_m(p))\neq (\mu_1(p), \dots , \mu_m(p))$. Without loss of generality, let $\lambda_1(p)\neq \mu_1(p)$. Then we have

$\mathcal{D}(\lambda)(\pi_1\mid_P)(p)=\lambda_1(p)\neq \mu_1(p)=\mathcal{D}(\mu)(\pi_1\mid_P)(p)$,\\
hence $\mathcal{D}(\lambda)\neq \mathcal{R}(\mu)$.

In order to prove that $\mathcal{D}$ is full, we first remark that $DMV_n$ is generated by the set of coordinate projections $\{\pi_1, \ldots \pi_n\}$.

Since $DMV_m \rightarrow \quot{DMV_m}{I(Q)}$ is a surjective homomorphism and $\quot{DMV_m}{I(Q)} \simeq DMV_m\mid_Q$, it is straightforward that  $\mathcal{D}(Q)=DMV_m\mid_Q$ is generated by $\{ \pi_1\mid_Q, \dots , \pi_m\mid_Q \}$.

Consider now $h: \mathcal{D}(Q)\rightarrow \mathcal{D}(P)$, with $P\subseteq [0,1]^n$, $Q\subseteq [0,1]^m$. We set

$\lambda_i=h(\pi_i\mid_Q)\in \mathcal{D}(P)$, $i=1, \dots , m$;\\
then, we define $\lambda$ as follows

$\lambda (p)=(\lambda_1(p), \dots \lambda_m(p))\subseteq [0,1]^m$, for any $p\in P$.\\
Trivially $\lambda$ is a $\mathbb{Q}$-map, since any $\lambda_i \in DMV_n\mid_{P}$, that is $\lambda_i$ is actually the restriction to $P$ of an element of $DMV_n$. We need to show that $\lambda (P)\subseteq Q$.

By Theorem \ref{thm:zerosetsDMV} there exists $f\in DMV_m$ such that $Q=V(f)$. Therefore $\lambda(P)\subseteq Q$ if $f(\lambda(P))=0$.\\
Since $f\in DMV_m$, and $DMV_m$ is freely generated by $\{\pi_1, \ldots \pi_m\}$, there exists a DMV-term $\sigma(\pi_1, \dots , \pi_m)$ such that $f=\sigma$. Moreover, $f\mid_Q=0$ implies \mbox{$h(f\mid_Q)=0$}. Hence
\vspace*{-0.2cm}
\begin{center}
$0=h(f\mid_Q)=\sigma (h(\pi_1\mid_Q), \dots , h(\pi_m\mid_Q))= \sigma (\lambda_1, \dots \lambda_m)$,
\end{center}
\vspace*{-0.2cm}
then $f(\lambda(p))=0$ for any $p\in P$ and the claim is settled.
\end{proof}

The following lemma is \cite[Corollary 3.6]{LucaEnzo} specialized in our case. We omit the proof, since it is again very similar to the aforementioned result. 
\blem
For a DMV-algebra $A$, the following are equivalent.

i) $A$ is finitely presented and projective,

ii) $A$ is finitely generated and projective,

iii) if $P \subseteq [0,1]^n$ is a polyhedron such that $A=\D(P)$, then $P$ is a retract of $[0,1]^n$ by a unital $\mathbb{Q}$-map.
\elem

\bfact
Trivially, the subcategory of  $\mathbf{RatPol_{[0,1]}^{\mathbb Q}}$ whose maps are $\mathbb{Z}$-maps is equivalent to the category  $\mathbf{RatPol_{[0,1]}^{\mathbb Z}}$, which is dual to the full subcategory of finitely presented MV-algebras.
\efact

Our last step is to extend the duality in Theorem \ref{teo:BBfpDMV} to rational polyhedra laying in $\mathbb{R}^n$, for some $n\in \mathbb{N}$. In order to do so, let us denote by $\mathbf{RatPol_{\mathbb{R}}^{\mathbb Q}}$ is the category of rational polyhedra laying in some appropriate $\mathbb{R}^n$ with $\mathbb{Q}$-maps, that is maps $\lambda=(\lambda_1, \dots , \lambda_m):P \subseteq \mathbb{R}^n\to Q\subseteq \mathbb{R}^m$ in which the components $\lambda_i: \mathbb{R}^n \to \mathbb{R}$ are piecewise linear maps with rational coefficients. In analogy with the unital case, we denote the set of continuous piecewise linear functions $f:\mathbb{R}^n \to \mathbb{R}$ by $DG_n$.

We then define a functor $\mathcal{P} :\mathbf{RatPol_{[0,1]}^{\mathbb Q}}\rightarrow \mathbf{RatPol_{\mathbb{R}}^{\mathbb Q}}$ by
\begin{itemize}
\item[]\textit{Objects:} For any $P\subseteq [0,1]^n$, $\mathcal{P} (P)=P$;
\item[]\textit{Arrows:} For any $\lambda: P\subseteq [0,1]^n \rightarrow Q\subseteq [0,1]^m$, we define $\mathcal{P} (\lambda)=\lambda$.
\end{itemize}

\blem \label{lem:P}
 $\mathcal{P}$ is a well defined functor.
\elem
\begin{proof}
We start by recalling that any piecewise linear function can be written as sups of infs of linear polynomials \cite[Theorem 2.1]{Ovc}. In particular, such linear functions are its components. That is, any piecewise linear function with either rational or real coefficients defined on some $[0,1]^n$ can be extended to $\mathbb{R}^n$ and $DMV_n=DG_n^r=\{ f\mid_{[0,1]^n} \mid f\in DG_n \}$. Hence, any $\lambda: P\subseteq [0,1]^n \rightarrow Q\subseteq [0,1]^m$ is indeed a morphism in $\mathbf{RatPol_{\mathbb{R}}^{\mathbb Q}}$, since we can find extensions $\overline{\lambda_i}:\mathbb{R}^n \to \mathbb{R}$ of any $\lambda_i$, $i=1, \dots, m$, such that $\sigma$ is the restriction to $P$ and the co-restriction to $Q$ of $(\overline{\lambda_1}, \dots, \overline{\lambda_m})$.

Finally, it is an easy exercise to show that $\mathcal{P}$ is indeed a functor.
\end{proof}

\bthm \label{theq1}
The functor $\mathcal{P}$ establishes a categorical equivalence between $\mathbf{RatPol_{[0,1]}^{\mathbb Q}}$ and $\mathbf{RatPol_{\mathbb{R}}^{\mathbb Q}}$.
\ethm
\begin{proof}
It is enough to prove that $\mathcal{P}$ is full, faithful and essentially surjective.

By \cite[Claim 3.5]{LucaEnzo}, for any rational polyhedron $P\subseteq \mathbb{R}^n$ there exists $d\ge 0$, a polyhedron $Q\subseteq [0,1]^d$ and a $\mathbb{Z}$-homeomorphism $\lambda:P\to Q$, hence $P\simeq Q=\mathcal{P}(Q)$ and the functor is essentially surjective. 

Fullness and faithfulness are straightforward by definition of $\mathcal{P}$, and the same arguments as in the proof of Lemma \ref{lem:P}.
\end{proof}

\bcor
$\mathbf{RatPol_{\mathbb{R}}^{\mathbb Q}}$ and $\mathbf{DMV_{fp}}$ are dual categories.
\ecor
\begin{proof}
It follows from Theorems \ref{teo:BBfpDMV} and  \ref{theq1}.
\end{proof}

The following figure summarizes the results related to geometrical dualities for DMV-algebras. For the categories of algebras, $\mathbf{fgp}$ stands for finitely generated and projective, $\mathbf{fp}$ stands for finitely presented, $\mathbf{fpp}$ stands for finitely presented and projective. $\mathbf{rRatPol^{\mathbb{Z}}_{[0,1]}}$ and $\mathbf{rRatPol^{\mathbb{Q}}_{[0,1]}}$ are the categories of polyhedra in some unit cube that are retract of the cube, with $\mathbb{Z}$-maps and $\mathbb{Q}$-maps respectively. $\mathbf{RatPol^{\mathbb{Z}}_{[0,1]}}$ and $\mathbf{RatPol^{\mathbb{Q}}_{[0,1]}}$ are the categories of polyhedra in some unit cube, with $\mathbb{Z}$-maps and $\mathbb{Q}$-maps respectively. $\mathbf{RatPol^{\mathbb{Z}}_{\mathbb{R}}}$ and $\mathbf{RatPol^{\mathbb{Q}}_{\mathbb{R}}}$ are the categories of polyhedra in some $\mathbb{R}^n$, with $\mathbb{Z}$-maps and $\mathbb{Q}$-maps respectively.\\

\begin{center}
\begin{tikzpicture}
  \node (rRatPolunitZ) {{$\mathbf{rRatPol^{\mathbb{Z}}_{[0,1]}}$}};
  \node (RatPolunitZ) [right=1cm of rRatPolunitZ] {{$\mathbf{RatPol_{[0,1]}^{\mathbb{Z}}}$}};
  \node (RatPolZ) [right=1cm of RatPolunitZ] {{$\mathbf{RatPol_{\mathbb{R}}^{\mathbb{Z}}}$}};
  \node (rRatPolunitQ) [below=1.5cm of rRatPolunitZ]{{$\mathbf{rRatPol^{\mathbb{Q}}_{[0,1]}}$}};
  \node (RatPolunitQ) [below=1.5cm of RatPolunitZ] {{$\mathbf{RatPol_{[0,1]}^{\mathbb{Q}}}$}};
  \node (RatPolQ) [below=1.5cm of RatPolZ] {{$\mathbf{RatPol_{\mathbb{R}}^{\mathbb{Q}}}$}};
  \draw [right hook->] (rRatPolunitZ) to node {}(rRatPolunitQ);
  \draw [right hook->] (RatPolunitZ) to node {}(RatPolunitQ);
  \draw [right hook->] (RatPolZ) to node {}(RatPolQ);
  \draw [right hook->] (rRatPolunitZ) to node {}(RatPolunitZ);
  \draw [right hook->] (RatPolunitZ) to node {}(RatPolZ);
  \draw [right hook->] (rRatPolunitQ) to node {}(RatPolunitQ);
  \draw [right hook->] (RatPolunitQ) to node {}(RatPolQ);
  
  \node (DMVfpg) [below=1.5cm of rRatPolunitQ]{{$\mathbf{DMV_{fgp}}$}}; 
  \node (DMVfpp) [below=1.5cm of RatPolunitQ]{{$\mathbf{DMV_{fpp}}$}};
  \node (DMVfp) [below=1.5cm of RatPolQ]{{$\mathbf{DMV_{fp}}$}};
  
  \draw [<->] (DMVfpg) to node {}(DMVfpp);
  \draw [right hook->] (DMVfpp) to node {}(DMVfp);
  \draw [<->] (DMVfpp) to node {}(rRatPolunitQ);
  \draw [<->] (DMVfp) to node {} (RatPolunitQ);
  \draw [<->] (DMVfp) to node {} (RatPolQ);
 
  \node (A) [below=0.3 of DMVfpp] {Figure 1. DMV-algebras and rational polyhedra};
\end{tikzpicture}
\end{center}

\section*{Acknowledgment}
S. Lapenta acknowledges partial support from the Italian National Research Project (PRIN2010-11) entitled \textit{Metodi logici per il trattamento dell’informazione}. I. Leu\c stean was supported by a grant of the Romanian National Authority for Scientific Research and Innovation, CNCS-UEFISCDI, project number PN-II-RU-TE-2014-4-0730.

\end{document}